\documentclass[graybox]{svmult}

\usepackage{type1cm}        
\usepackage{makeidx}         
\usepackage{graphicx}        
\usepackage{multicol}        
\usepackage[bottom]{footmisc}

\usepackage{newtxtext}       %
\usepackage[varvw]{newtxmath}       

\usepackage{cleveref}
\usepackage{booktabs}

\makeindex             

\newcommand{\bbC}{\mathbb{C}}
\newcommand{\bbR}{\mathbb{R}}
\newcommand{\bbN}{\mathbb{N}}

\newcommand{\scrM}{\mathscr{M}}

\newcommand{\mcal}[1]{\mathcal{#1}}

\newcommand{\rank}{\mathrm{rank}}
\newcommand{\cprank}{\rank_{\text{cp}}}
\newcommand{\seprank}{\rank_{\text{sep}}}

\newcommand{\tcp}{\tau_{\text{cp}}}

\newcommand{\xicp}{\xi^{\mathrm{cp}}}
\newcommand{\xisp}{\xi^{\mathrm{isp}}}
\newcommand{\xicpsp}{\xi^{\mathrm{cp,isp}}}
\newcommand{\xicpwsp}{\xi^{\mathrm{cp,wisp}}}

\newcommand{\conv}{\mathrm{conv}}
\newcommand{\supp}{\mathrm{supp}}

\newcommand{\val}{\mathrm{val}}
\newcommand{\valpop}{\mathrm{val}^\mathrm{POP}}
\newcommand{\valsp}{\mathrm{val}^\mathrm{sp}}

\newcommand{\wtl}{\widetilde} 
\newcommand{\wh}{\widehat}
\newcommand{\meas}{\scrM}

\newcommand{\inp}[1]{\langle #1 \rangle}

\begin{document}

\title*{Matrix factorization ranks via polynomial optimization}
\author{Andries Steenkamp
\thanks{Centrum Wiskunde \& Informatica (CWI), Amsterdam. \url{andries.steenkamp@cwi.nl}
\newline
This work is supported by the European Union's Framework Programme for Research and Innovation Horizon
2020 under the Marie Skłodowska-Curie Actions Grant Agreement No. 813211  (POEMA).
}}
\institute{Andries Steenkamp \at CWI Amsterdam, \email{jajs@cwi.nl}}
%
%
\maketitle

\abstract{
In light of recent data science trends, new interest has fallen in alternative matrix factorizations. By this, we mean various ways of factorizing particular data matrices so that the factors have special properties and reveal insights into the original data. We are interested in the specialized ranks associated with these factorizations, but they are usually difficult to compute. In particular, we consider the nonnegative-, completely positive-, and separable ranks. We focus on a general tool for approximating factorization ranks, the moment hierarchy, a classical technique from polynomial optimization, further augmented by exploiting ideal-sparsity. Contrary to other examples of sparsity, the resulting sparse hierarchy yields equally strong, if not superior, bounds while potentially delivering a speed-up in computation.
}

\section{Introduction and motivation for matrix factorization ranks}
\label{sec:1}
We live in a digital world. Data drives decisions as a never-ending stream of information engulfs our lives. An
essential tool for navigating the flood of information is the ability to distill large bodies of information
into actionable knowledge. A practical example is \emph{nonnegative (NN) factorization}, applied to data represented as a matrix. A NN factorization of an entry-wise nonnegative matrix $M \in \bbR^{m \times n}_+$ is a pair of nonnegative matrices $A\in \bbR^{m \times r}_+$ and $B\in \bbR^{r \times n}_+$ for some integer $r\in \bbN$ such that:
\begin{equation} \label{MAB}
M = AB.
\end{equation}
The primary object of interest is the inner dimension  $r$ of the factorization. One can always take $A=M$ and $B=I$, where $I$ is the identity matrix, hence getting $r=n$. However, the interesting case is when $r < \frac{mn}{m+n}$. In this case, one has managed to express the $m\times n$ values of $M$ in terms of the $(m \times r) + (r \times n)$ values of $A$ and $B$, and as a result, using less storage. The smallest integer $r$ for which this is possible is called the \emph{nonnegative matrix factorization rank}, or just the \emph{nonnegative rank} for short, and is mathematically defined as follows,
\begin{equation} \label{nnrankdef}
\mathrm{rank}_+(M) := \min \{r \in \mathbb{N} : M = AB \text{ for some } A\in \bbR^{m \times r}_+ \text{ and } B\in \bbR^{r \times n}_+  \}.
\end{equation}
It is not hard to see that the NN-rank is sandwiched between the classical rank and the size of the matrix, i.e.,
$$
\rank(M) \leq \rank_+(M) \leq \min\{n,m\}.
$$
However, storage space efficiency is only part of the value of NN factorization. The true power of NN factorization comes from the fact that it is an easy-to-interpret \emph{linear dimensionality reduction} technique. To understand what we mean by this, we first re-examine the relationship between the three matrices $M, A$, and $B$ in \cref{MAB}.

Observe how the $j^{\mathrm{th}}$ column of $M$ is given as a conic combination of the columns of $A$ with weights given by the  $j^{\mathrm{th}}$ column of $B$, i.e.,
\begin{equation} \label{nnmf_interp}
M_{:,j} = \sum_{i=1}^r B_{i,j}A_{:,i}.
\end{equation}
Because all terms involved are nonnegative, zero entries in $M$ force the corresponding entries of the factors to be zero.
Formally, for any $k\in \{1,2,..,m\}$ and $j \in \{1,2,..,n\}$, $M_{k,j} = 0$ if and only if $B_{i,j}A_{k,i} = 0$ for all $i \in \{1,2,..,r\}$. 
Having no cancellation among factors will be useful for interpreting applications of nonnegative factorization. We will explain this further in \cref{sec:1:2} with examples. Furthermore, observe that the nonnegative factorization needs not be unique. In fact, for any non-singular, nonnegative matrix $P \in \bbR^{r \times r}_+$ with nonnegative inverse $P^{-1}$ one can produce another factorization
\begin{equation} 
M = (AP^{-1})(PB).
\end{equation}
An example of such a matrix $P$ would be a permutation matrix.
\begin{example}
{Example of a nonnegative factorization}
Consider the following example of a $4\times 4$ nonnegative matrix and its nonnegative factorization
from which we can deduce that $\rank_+(M)=2$, because $ \rank(M)=2$:
{
$$ M =   
  \left[ {\begin{array}{cccc}
    35  & 38  & 41 &  44\\
    79  & 86  & 93 & 100\\
    123 & 134 & 145&  156\\
    167 & 182 & 197&  212\\
  \end{array} } \right]
  = 
  \left[ {\begin{array}{cc}
  1& 2\\
  3& 4\\
  5& 6\\
  7& 8\\ 
  \end{array} } \right]
  \left[ {\begin{array}{cccc}
  	 9  &10 &11 &12 \\
     13 &14 &15 &16 \\ 
  \end{array} } \right] = AB.
$$}
\end{example}

\subsection{Applications of nonnegative factorization}
\label{sec:1:2} 
Having introduced the nonnegative rank, we now justify its importance with three applications. What we present here is but a small fraction of the whole body of literature on nonnegative factorization. The interested reader is highly encouraged to read a recent monograph of Gillis \cite{doi:10.1137/1.9781611976410} for an in-depth study of the nonnegative rank with many applications and further references. Alternatively, we invite the reader to try and conceive a few applications of their own.

\paragraph{\textbf{Image processing}}
When analyzing large amounts of images, it is natural to ask if the vast bulk of images are not just combinations of a few ``basic images". This raises two questions. First, how does one find or construct a set of basic images? Second, given this, hopefully small, set of basic images, how does one reproduce the original images? Lee and Seung answered both questions in \cite{Lee1999LearningTP}, where they factorized a set of images of human faces into a set of typical facial features and nonnegative weights. Combining the weights and features, one approximately recovers the original faces. 
In this setting, the matrix $M$ has as columns the vectorized gray-scale images of human faces, hence $M_{i,j}$ is the $i^{\mathrm{th}}$ pixel of the $j^{\mathrm{th}}$ face, with a value between $0$ and $1$, with $0$ corresponding to black and $1$ to white.

Recalling the interpretation of (\ref{nnmf_interp}), we can think of the columns of the matrix $A$ as (vectorized) images of human facial features, like a mouth or pair of eyes.
Hence, the $j^{\mathrm{th}}$ image $M_{:,j}$ is a weighted sum of feature-images $A_{:,i}~(i \in [r])$, 
where the (nonnegative) weight of feature $A_{:,i}$ is given by entry $B_{i,j}$ of the matrix $B$.

 In contradistinction to other techniques like \emph{principle component analysis} (PCA), which possibly gives factors with negative entries, NN factorization saves us from the task of interpreting notions like ``negative pixels" or ``image cancellations." By ``negative pixels," we mean negative factor values, i.e., $B_{i,j}A_{k,i} < 0$. This means that factor $B_{i,j}A_{:,i}$ does not just add features but also possibly erases the features added by other factors $B_{\ell,j}A_{:,\ell}$, where $\ell \neq i$. A fun by-product of these image factorizations is that one can generate new images by multiplying the matrix $A$ with new weights different from $B$. However, the resulting images are not guaranteed to look like faces for a poor choice of weights.

\paragraph{\textbf{Topic recovery and document classification}}
In text analysis, the matrix $M$ is called the \emph{word occurrence matrix}, and its entries $M_{i,j}$ are the number of times the $i^{\mathrm{th}}$ word occurs in the  $j^{\mathrm{th}}$ document. This way of looking at a corpus of text is often called a ``bag of words model,"  the sequence is ignored, and only the quantity is considered. Since word count is always nonnegative, $M$ is nonnegative and has some NN factors $A$ and $B$. The columns of matrix $A$ take the meaning of ``topics", and $B$ gives the correct weights to recover $M$. Since there are no cancellations, we observe in the columns of $A$ that certain sets of words tend to occur together, at least within the original set of documents. Moreover, we see how the documents (the columns of $M$) are composed from these base topics (the columns of $A$), with the importance of each topic given by the entries of $B$. One can hence use these learned topics to group or classify documents.
 
\paragraph{\textbf{Linear extension complexity}}
This third application is different from the above two. First, we define the linear extension complexity, then we show how it relates to the nonnegative rank, and finally, we motivate its importance. The \emph{linear extension complexity} of a polytope $P$ is the smallest integer $r$ for which $P$ can be expressed as the linear image of an affine section of $\bbR_+^r$. Alternatively, the linear extension complexity can be defined as the smallest number of facets a higher dimensional polytope $Q$ can have while still having $P$ as a projection. In 1991 Yannakakis \cite{Yannakakis1988ExpressingCO} proved that the linear extension complexity of $P$ is equal to the nonnegative rank of the \emph{slack matrix}  associated with $P$. For a polytope $P$ the slack matrix is $(d_i-c_i^Tv)_{v \in \mcal{V}, i\in \mcal{I}}$, where $c_i\in\bbR^m,~d_i\in\bbR$ come from the hyperplane representation of $P = \{x \in \bbR^m: c^T_i x \leq d_i ~~(i \in \mcal{I}) \}$, and the vectors $v \in \bbR^m $ come from the extremal point representation of $P = \conv{\mcal(V)}$. This link between nonnegative rank and linear extension complexity was instrumental in showing why many combinatorial problems, like the traveling salesman problem, could not be efficiently solved simply by lifting the associated problem polytope to higher dimensions in some clever way, see \cite{10.1145/2716307}. On the topic of lifting convex sets we refer the reader to the survey \cite{doi:10.1137/20M1324417}.

\subsection{Commonly used notation}
We group here some notation used throughout the chapter. For any positive integer $m \in \bbN$ we denote by $[m] := \{1,2,...,m\}$ the set consisting of the first $m$ positive integers. For vectors $u,v \in \bbR^m$ we denote by $\inp{u,v} := \sum_{i\in[m]} u_iv_i$ the vector inner-product of $u$ and $v$. Similarly, for matrices $A, B \in \bbR^{m \times n}$ we use the same notation to denote the Frobenius inner product $\inp{A, B}:= \sum_{i\in[m], ~j\in[n]} A_{i,j} B_{i,j}$  of $A$ and $B$. We say that a square matrix $A \in \bbR^{m \times m}$ is \emph{positive semi-definite} (PSD), denoted by $A \succeq 0 $, if and only if $v^TAv \geq 0$ for any choice of $v \in \bbR^m$. The set of all PSD matrices of size $r \in \bbN$ is denoted by $\mcal{S}^{r}_+ := \{A \in \bbR^{r \times r} : A \succeq 0 \}$. Analogous to real matrices that are PSD, there are complex square matrices that are \emph{Hermitian PSD}. A complex matrix $A \in \bbC^{m \times m}$ is Hermitian PSD if and only if $v^*Av \geq 0$ for all  $v \in \bbC^m$, where $v^*$ is the complex conjugate of $v$. For positive integer $m$, we denote by $\mcal{H}^m$ the set of all $m\times m$ Hermitian matrices. By 
$$
\Sigma[x]:= \{ \sum_{i \in [k]}\sigma_i: k \in \bbN,~ \sigma_i = p_i^2 \text{ for some polynomial } p_i \in \bbR[x]\}
$$ 
we denote the set of all sums of squares of polynomials. If the variables $x=(x_1,x_2,...,x_m)$ are clear from the context we simply write $\Sigma$.

\subsection{On computing the nonnegative rank}
Given the utility of computing nonnegative factorizations, it is natural to ask the following. Is it difficult to compute the nonnegative rank for a given data matrix $M \geq 0$? This was answered in the affirmative in 2009 by Vavasis \cite{Vavasis2009OnTC}. 
Despite being NP-hard to solve, good approximations are sometimes quite accessible. In \cref{sec:2}, we show a general technique for approximating the matrix factorization rank from below using tools from polynomial optimization. An alternative geometrically motivated approach is to look for a minimal \emph{rectangle cover} for the support of $M$, see \cite{GILLIS20122685}.Given a matrix $M \in \bbR^{n \times m}$, one seeks the smallest set of \emph{rectangles}, sets of the form $R := \{ \{i,j\}: i \in I \subset [n],  j \in J \subset [n]\}$, such that for each nonzero entry $M_{i,j} \neq 0$, $\{i,j\}$ belongs to at least one of these rectangle. 

Finding the factorization rank does not necessarily give a factorization. The method we propose in \cref{sec:2} does not generally give a factorization, except in a particular case, which we consider in \cref{sec:2.3}, where it is possible to recover the factors. For NN factorization, several algorithms exist that iteratively compute $A$ and $B$ given a guessed value $r$. However, these algorithms only give approximate factorizations, that is, $M \approx AB$, with respect to some norm. A sufficiently good approximate NN factorization also implies an upper bound on the NN rank. 

 For practical problems, an approximation is often sufficient. For a detailed account of NN factorization, we refer the reader again to the book of Gillis \cite{doi:10.1137/1.9781611976410}.

\subsection{Other factorization ranks}
Above, we looked at NN factorization and some of its applications in data analysis and optimization theory. However, there are many more matrix factorization ranks, each having its intricacies, applications, and interpretations. We list a few more examples of factorization ranks to link them to NN factorization and later state some results on some of them.

\paragraph{\textbf{Completely positive factorization}}
This factorization is very similar to nonnegative factorization apart from the modification that $B=A^T$. Formally, an entry-wise nonnegative matrix $M \in \mathcal{S}^{m}$ is \emph{completely positive} (CP) if there exists a nonnegative matrix  $A\in \bbR^{m \times r}_+$, for some integer $r\in \bbN$, such that:
\begin{equation}
M = AA^T.
\end{equation}
Clearly, CP-matrices are \emph{doubly nonnegative}, i.e., entry-wise nonnegative and \emph{positive semi-definite} (PSD). However, these are not sufficient criteria for $M$ to be CP, unless $m \leq 4$, see \cite{doi:10.1142/5273}. Consider the following (Example 2.9 from  \cite{doi:10.1142/5273}) to see that this does not hold for $m\geq 5$. 

\begin{example}{Example of a doubly nonnegative matrix that is not CP \cite{doi:10.1142/5273}} 
\[
M =
\left[ {\begin{array}{ccccc}
    1& 1& 0& 0& 1\\
	1& 2& 1& 0& 0\\
	0& 1& 2& 1& 0\\
	0& 0& 1& 2& 1\\
	1& 0& 0& 1& 6\\
  \end{array} } \right].
 \]
The nonnegativity is clear, the PSDness is checked via computing all the minors (or using a computer to check the eigenvalues).
To see why the matrix is not CP, we refer to the explanation given in Berman and Shaked-Monderer's monograph \cite{doi:10.1142/5273}.
\end{example}

Just deciding if a given matrix is CP is already an NP-hard problem; see \cite{Dickinson2014OnTC}. Because the \emph{completely positive factors} $A$ and $A^T$ are the same, up to transposition, the CP factorization is often called a \emph{symmetric factorization}. We will see another example shortly at the end of this section. Similar to nonnegative rank, there is a \emph{completely positive rank} mathematically defined as the smallest inner dimension $r\in \bbN$ for which a CP factorization of $M$ exists, i.e.
\begin{equation} \label{cprankdef}
\cprank (M) := \min \{r \in \mathbb{N} : M = AA^T \text{ for some } A\in \bbR^{m \times r}_+ \}.
\end{equation}
Clearly a matrix $M$ is CP if and only if $\cprank (M) < \infty$. Hence computing the CP-rank can't be any easier than deciding if $M$ is CP. That being said, the complexity status of computing $\cprank(M)$ for a CP matrix $M$ is unknown to the best of our knowledge. Some upper bounds are known for the CP-rank: first, $\cprank(M)  \leq m$,  when $m \leq 4$, and second, $\cprank(M) \leq {m+1 \choose 2} -4 $ if $m\geq 5$, see \cite{ShakedMonderer2013OnTC}. In 1994 it was conjectured by Drew, Johnson, and Loewy \cite{Drew1994CompletelyPM} that $\cprank(M) \leq  \lfloor \frac{m^2}{4} \rfloor$, the bound being only attained for CP matrices $M$ that have complete bipartite support graphs. This conjecture was disproved by Bomze et al. \cite{BOMZE2014208,doi:10.1137/140973207} two decades later, using several specially constructed counter-examples. We show in \cref{M7} an example, namely $\wtl M_7$ from \cite{BOMZE2014208}, of size $m=7$, with $\cprank (\wtl M_7)= 14 > \lfloor \frac{49}{4} \rfloor = 12$:
\begin{equation} \label{M7} 
  \wtl M_7 =
  \left[ {\begin{array}{ccccccc}
    163& 108& 27& 4 &4& 27& 108 \\
	108& 163& 108& 27 &4 &4& 27 \\
	27& 108& 163 &108& 27& 4& 4 \\
	4& 27 &108 &163 &108& 27& 4 \\
	4& 4 &27 &108& 163& 108& 27 \\
	27 &4 & 4& 27 &108& 163& 108 \\
	108& 27& 4 &4 &27 &108& 163 
  \end{array} } \right].
\end{equation} 

On the more applied side, CP matrices occur in the theory of \emph{block designs}. We omit many details here, but essentially, block designs deal with arranging distinct objects into blocks in such a way that the objects occur with certain regularity within and among the blocks. There is a direct application of block design in designing experiments where researchers wish to prevent the differences between test subjects from obfuscating the differences in outcome due to treatment, see \cite{Hall1988CombinatorialTH} for block designs in depth and see \cite{ShakedMondererCopositiveAC} for the link between block designs and CP matrices. 
 
From another perspective of applications, completely positive matrices are of great interest in optimization.

In 2009, Burer \cite{burerCopositiveRepresentationBinary2009} showed that any nonconvex quadratic program with binary and continuous variables could be reformulated as a linear program over the cone of completely positive matrices. This effectively meant that many NP-hard problems could now be viewed as linear programs with CP-membership constraints. This reformulation does not make the problems any easier to solve as the difficulty is now pushed into characterizing complete positivity. However, it does allow us to attack a large class of problems by understanding a unifying thread. 

For a thorough account of completely positive and copositive matrices, we refer the inquisitive reader to the monograph by Berman and Shaked-Monderer \cite{ShakedMondererCopositiveAC}.

\paragraph{\textbf{Separable rank}}
In the setting of quantum information theory, the state of some physical system is often characterized by a Hermitian PSD matrix $M\in \mcal{H}^m \otimes \mcal{H}^m$. These states are said to be \emph{separable} if there exist vectors $a_1,...,a_r,b_1,...,b_r \in \bbC^m$ for which 
\begin{equation}
M = \sum_{\ell=1}^r a_\ell a_\ell^* \otimes b_\ell b_\ell^*,
\end{equation}
where $a^*$ denotes the complex conjugate of $a$, and $\otimes$ denotes the tensor product. We will not go into the details, but it suffices to think of separable states as fully explained by classical physics, in contradistinction, non-separable states, a.k.a.  \emph{entangled states} have special properties of interest in quantum physics. For rank-one states, i.e., if $\rank(M)=1$, also called \emph{pure states}, one can obtain a separable factorization by using \emph{singular value decomposition (SVD)}. Non-rank-one states are called \emph{mixed states}, and deciding whether a mixed state $M$ is separable is, in general, NP-hard, see \cite{10.1145/780542.780545,Gharibian2010StrongNO}. 

\begin{example}{Example of an entangled state \cite{choiPositiveLinearMaps1982}} 
Consider the following mixed state of size $9 \times 9$, hence $m=3$. We have omitted showing zeros for readability, and we draw grid lines in order to highlight the block structure.
$$
M =
	\left[
	\begin{array}{ccc|ccc|ccc}
		1& & & &1& & & &1\\
		&2& &1& & & & & \\
		& &\frac{1}{2}& & & &1& & \\
		\hline
		&1& &\frac{1}{2}& & & & & \\
		1& & & &1& & & &1\\
		& & & & &2& &1& \\
		\hline
		& &1& & & &2& & \\
		& & & & &1& &\frac{1}{2}& \\
		1& & & &1& & & &1
	\end{array}
	\right].
	$$
In \cite{choiPositiveLinearMaps1982} it is shown that $M$ is entangled.	
	
\end{example}

Analogously to other matrix ranks we considered thus far, there is also a notion of \emph{separable rank}, see \cite{DelasCuevas2019separabilitymixed}, sometimes called \emph{optimal ensemble cardinality} \cite{doi:10.1080/09500340008244049}, which we define for a separable matrix $M$ as 
\begin{equation} \label{sperankdef}
\seprank(M) = \min \{r \in \bbN : M = \sum_{\ell=1}^r a_\ell a_\ell^* \otimes b_\ell b_\ell^* \text{ for some } a_\ell \text{ and} \ b_\ell  \text{ in } \bbC^{m} \}.
\end{equation}
A possible interpretation of the separable rank is that it gives a sense of how complex a classical system is, with the convention being that an entangled state has infinite separable rank. To our knowledge, the complexity of computing the separable rank is still unknown. There are some crude bounds on the separable rank though
$$
\rank(M) \leq \seprank(M) \leq \rank(M)^2 .
$$
The left most inequality can be strict (see \cite{doi:10.1080/09500340008244049}) and the right most inequality follows from Caratheodory’s theorem \cite{uhlmannEntropyOptimalDecompositions1998}.

In addition to the above definition, there are several other variations on this notion of separability. One variation is to look for factorizations of the form $M = \sum_{\ell=1}^r A_\ell \otimes B_\ell $, where $A_\ell, B_\ell\in \mcal{H}^m$ are Hermitian PSD matrices. From this it is easy to define the associated \emph{mixed separable rank} as the smallest $r$ for which such a factorization is possible. When $M$ is diagonal its mixed separable rank equals the nonnegative rank of an associated $m \times m$ matrix consisting of the diagonal entries of $M$, see \cite{delasCuevas2020MixedSI}. This shows that mixed separable rank is hard to compute.

\paragraph{\textbf{Nonnegative tensor factorization ranks}} 
Tensors, also called multi-way arrays, are natural generalizations of matrices commonly encountered in applied fields such as statistics, computer vision, and data science. Strictly speaking, tensor factorization ranks falls beyond the scope of this chapter, but given the similarities, we would be remiss not to include some remarks and references on the matter.

Consider, for example, a three-way array $T \in \bbR^{n\times m\times p}$. Then, its \emph{tensor rank} is the smallest number $r\in \bbN$ of rank-one tensors (tensors of the form $a\otimes b \otimes c$ for some $a \in \bbR^n$, $b \in \bbR^m$, and $c \in \bbR^p$) necessary to describe $T$, i.e.
$$
\rank_{\mathrm{tensor}}(T) = \min \{r \in \bbN  : T = \sum_{\ell=1}^r a_\ell \otimes b_\ell \otimes c_\ell  ,~a_\ell \in \bbR^n,~ b_\ell \in \bbR^m,~ c_\ell \in \bbR^p \}.
$$ 
Similarly, one can define the \emph{nonnegative tensor rank} by requiring the factors $a_\ell,b_\ell,$ and $c_\ell$ to be nonnegative. Moreover, one can define the \emph{symmetric tensor rank} by requiring $n=m=p$ and the factors to be equal, i.e., $a_\ell=b_\ell=c_\ell$ for all $\ell \in [r]$. 
An interesting effect of going to tensors is that some decompositions become unique \cite{https://doi.org/10.1002/1099-128X(200005/06)14:3<229::AID-CEM587>3.0.CO;2-N}. See \cite{10.5555/1822971} for an applications-centric monograph on tensor factorization. For a mathematical survey, see Kolda and Bader \cite{doi:10.1137/07070111X}.

\paragraph{\textbf{Non-commutative matrix factorization ranks}} 
We conclude this section with two more factorizations, often called non-commutative analogs of nonnegative- and CP factorizations.
First is the \emph{positive semi-definite (PSD) factorization}, where given a matrix $M\in \bbR^{m \times n}_+$ we look for PSD matrices $A_1,...,A_m,B_1,...,B_n \in \mcal{S}^{r}_+$, for some $r \in \bbN$, such that the matrix $M$ is described entry-wise as follows: $M_{i,j} = \inp{A_i,B_j} \in \bbR$ for $i\in [m]$ and $j\in [n]$. If the matrices $A_i$ and $B_j$ are diagonal for $i \in [m]$ and $j \in [n]$, then we recover a nonnegative factorization. Similar to nonnegative factorization, there is a substantial research interest in PSD-factorization, largely due to its many appealing geometric interpretations, including semi-definite representations of polyhedra. We refer the reader to the survey by Fawzi et al. \cite{fawziPositiveSemidefiniteRank2015} for further study of PSD-factorizations.

Second, we have the symmetric analog of PSD-factorization, called a \emph{completely positive semi-definite (CPSD) factorization}, which simply adds the requirement that $n=m$ and $B_i =A_i$ for all $i \in [m]$. 

The associated definition of rank accompanying these two factorizations should be clear. 
We stop introducing matrix factorization ranks now and shift gears towards proving bounds.

\section{Bounding matrix factorization ranks}
\label{sec:2}
There are two modes of approximating factorization ranks. The first is from above, using heuristics to construct factorizations. The second is from below, via computing parameters, often combinatorial in nature, exploiting the support graph of the matrix $M$. The approach we follow falls in this latter category. We will explain the method applied to the CP-rank, though it should be clear what substitutions are needed to generalize it to the other factorization ranks we described in \cref{sec:1}.

In this section, we give the focal point of this chapter, the \emph{moment hierarchy}. This is the core technique for approximating factorization ranks. In order to define the hierarchy, we start in \cref{sec:2:1} with polynomial optimization problems, which we feel is more natural than jumping straight into generalized moment problems considered in \cref{sec:2.2}. Finally, we apply the tools built in \cref{sec:2.2} to the setting of CP-rank in \cref{sec:2.3}.

The significant results concerning the properties of the hierarchy, like convergence, flatness, and exactness, will be mentioned and explained as we proceed. See the works \cite{Gribling2019LowerBO, Gribling2022BoundingTS} for a more fleshed-out exposition of the process we follow here.

\subsection{A brief introduction to polynomial optimization}
\label{sec:2:1}
We will be drawing heavily from the rich field of polynomial optimization, so it is only natural that we introduce some tools and notations in this regard. This small subsection is not an overview of the field. For that, we recommend the excellent works \cite{Laurent2009,doi:10.1142/p665}. We attempt to cover only the necessities needed to motivate the title and ease the reader into the more advanced machinery. Depending on the book's other chapters and the reader's background, some topics may be familiar, in which case, perusing this subsection will at least clarify the notation we use.

Consider the following optimization problem:
\begin{equation} \label{POP}
f^{\mathrm{min}} := \inf_{x \in K} f(x), 
\end{equation}
where 
\begin{equation} \label{semi_alg}
K := \{ x \in \bbR^{m} : g_i(x) \geq 0 ~(i \in [p]) ~,~h_j(x) = 0 ~(j \in [q]) \},
\end{equation}
and $f,g_1,...,g_p,h_1,...,h_q \in \bbR[x]$ are polynomials in $m$ variables $x_1,...,x_m$. The domain of optimization, $K$, is a basic closed semi-algebraic set. Problem (\ref{POP}) is called a \emph{polynomial optimization problem} or a \emph{POP} for short. POPs are versatile tools for modeling various problems. Clearly, linear and quadratic programs are instances of POPs. Furthermore, one can encode binary variables with polynomial constraints of the form $x_i(x_i-1)=0$. Hence, many NP-hard problems can be reformulated as a POP, and, as such, POPs are generally hard to solve, see \cite{doi:10.1142/p665}.

The moment approach to attacking a POP of the form (\ref{POP}) is as follows. We optimize the integral of the objective over probability measures that have support on $K$, i.e., we consider the following problem
\begin{equation}  \label{measPOP}
\begin{split} 
\valpop := &\inf_{\mu \in \scrM(K)} \int f(x) d\mu,  \\
        &~~~~~~\text{s.t.}~~ ~\int 1 d\mu =1,
\end{split} 
\end{equation}
where $\scrM(K)$ is the set of all Borel measures supported on the set $K$.

Problems (\ref{POP}) and (\ref{measPOP}) are equivalent in the sense that they have the same optimal values, i.e., $f^{\mathrm{min}} = \valpop$. To see that $f^{\mathrm{min}} \geq \valpop$ holds consider the Dirac delta measure $\delta_{x^{\mathrm{min}}}$ supported at a
minimizer $x^{\mathrm{min}}$ of problem (\ref{POP}). Then we have
$$
\valpop \leq  \int f(x) d\delta_{x^{\mathrm{min}}} =   f(x^{\mathrm{min}}) = f^{\mathrm{min}}.
$$
For the other inequality, $\valpop \geq f^{\mathrm{min}}$, simply observe
$$
\int f(x) d\mu  \geq f^{\mathrm{min}} \int 1 d\mu =  f^{\mathrm{min}}.
$$
The last equality comes from the fact that $\mu$ is a probability measure. With the equivalence between POPs and this new problem (\ref{measPOP}) established, we can focus on solving the latter.

\subsection{\textbf{Generalized moment problems}}
\label{sec:2.2}
Problem (\ref{measPOP}) is a special instance of what is called a \emph{generalized moment problem} (GMP), which is an even more versatile type of problem than a POP. Reformulating our optimization problem over measures does not give a clear advantage, as measures are difficult to handle. However, we will soon see in (\ref{hier}) how one can truncate the problem to create a hierarchy of \emph{semi-definite programs} (SDP).

Consider the following general form of GMP:
\begin{equation}  \label{GMP}
\val := \inf_{\mu \in \scrM(K)} \big{\{}  \int f_0(x) d\mu : \int f_i d\mu = a_i ~(i \in [N]) \big{\}},
\end{equation}
where$f_0,f_1,...,f_N$ are polynomials. From the discussion in \cref{sec:2:1}, we saw that POPs are a special class of GMPs with $N=1$ and $f_1 = a_1 = 1$. In \cref{sec:2.3}, we will show that the CP-rank can be reformulated as a GMP. 
Before we can start attacking the above GMP with the so-called \emph{moment method}, we must first set some notation and basic definitions. Let $\bbN^m_t$ be the set of all \emph{multi-indices} $\alpha = (\alpha_1,\alpha_2,...,\alpha_m)$ such that $|\alpha|:= \sum_{i =1}^m \alpha_i \leq t $. If $t = \infty$ we just write $\bbN^m$. For $x= (x_1,x_2,...,x_m) \in \bbR^m$ denote the truncated sequence of monomials by $[x]_t := \big{(} x^\alpha \big{)}_{\alpha \in \bbN^m_t}$, where $x^\alpha:= x_1^{\alpha_1}x_2^{\alpha_2} \cdots x_m^{\alpha_m}$. For a measure $\mu$ define its \emph{moments} to be the sequence of values, obtained when integrating the monomials w.r.t. the measure, i.e,
\begin{equation*} \label{defmoments}
\int x^{\alpha} d\mu ~~ (\alpha \in \bbN^m).
\end{equation*}
Using moments, we can think of measures as \emph{linear functionals} acting on the ring of polynomials. That is, for a measure $\mu$, we can define a linear map $L: \bbR [x] \to \bbR$ by defining what it does to monomials: $L(x^\alpha) = \int x^{\alpha} d\mu $ for every $\alpha \in \bbN^n$. Hence, for any polynomial $f = \sum_{\alpha} c_\alpha x^\alpha $,  we have
$$
L(f) = L(\sum_{\alpha} c_\alpha x^\alpha )= \sum_{\alpha} c_\alpha L(x^\alpha) = \sum_{\alpha} c_\alpha \int x^{\alpha} d\mu = \int f d\mu.
$$
Denote the space of truncated linear functionals acting on the space of polynomials of degree at most $t$, i.e., on $\bbR[x]_t$, by $\bbR[x]^*_t$.
Going in the opposite direction, i.e., from a linear functional $L \in \bbR[x]^*$ to a measure $\mu$, is not always possible.
When there does exists a measure $\mu \in \scrM(K)$ such that $L(x^\alpha) = \int x^{\alpha} d\mu $ for all $\alpha \in \bbN^m$, $L$ is said to have a \emph{representing measure}. We introduce some more concepts and notation to characterize the necessary conditions for measure representable functionals.

Recall the definition of our semi-algebraic set $K$ in \cref{semi_alg}:
\begin{equation*} 
K := \{ x \in \bbR^{m} : g_i(x) \geq 0 ~(i \in [p]) ~,~h_j(x) = 0 ~(j \in [q]) \}.
\end{equation*}

For $t \in \bbN \cup \{\infty\}$ define the \emph{truncated quadratic module} generated by ${ \bf g }:= (g_0,g_1,g_2,...,g_p)$, with $g_0=1$, as
\begin{equation} \label{quadmod}
\mcal{M}({ \bf g })_{2t} := \big{\{} \sum_{j = 0}^p \sigma_j g_j : \sigma_j \in \Sigma,~ \deg(\sigma_j g_j) \leq 2t \big{\}}. 
\end{equation}
Here $\Sigma$ denotes the set of \emph{sums of squares of polynomials}. In a similar vein we define the truncated \emph{ideal} generated by ${ \bf h }:= (h_1,h_2,...,h_q)$ as
\begin{equation} \label{ideal}
\mcal{I}({ \bf h })_{2t} := \big{\{} \sum_{j=1}^{q} \gamma_j h_j : \gamma_j \in \bbR[x],~ \deg(\gamma_j h_j) \leq 2t \big{\}}.
\end{equation}
When $t = \infty$ we also just drop the subscript and write: $\mcal{M}({ \bf g })$ and $\mcal{I}({ \bf h })$.

A crucial observation (see \cite{doi:10.1142/p665}) is that  if $L \in \bbR[x]^*$ has a representing measure $\mu \in \scrM(K)$ then $L \geq 0$ on $\mcal{M}({ \bf g })$ and $L = 0$ on $\mcal{I}({ \bf h })$. Using this, we can define, for any $t \in \bbN \cap \{\infty\}$, the following sequence of parameters,

\begin{equation}  \label{hier}
\begin{split} 
\xi_t:=  \min \{ L(f_0) :& L \in \bbR[x]^*_{2t},\\
			 & L(f_i) = a_i ~~(i \in [N]),\\
			 & L \geq 0 \text{ on } \mcal{M}({ \bf g })_{2t}, \\
			 & L = 0 \text{ on } \mcal{I}({ \bf h })_{2t} \}.
\end{split}
\end{equation}
We call $\xi_1,\xi_2,...,\xi_\infty$ a hierarchy as we clearly have for any $t \in \bbN$ that, 
$$
\xi_t \leq \xi_{t+1} \leq \xi_\infty \leq \val .
$$

\begin{example}{Exercise 1} 
Given a solution $L$ to the problem associated with $\xi_{t+1}$ construct a solution to the problem associated with $\xi_{t}$.
\end{example}

Under mild assumptions, the bounds $\xi_t$ converge asymptotically to
the optimum value $\val$ as $t$ goes to infinity. We state the well-known and widely used result here and refer to \cite{doi:10.1142/p665, Klerk2019ASO} for a full exposition.

\begin{theorem}\label{theoconvGMP}
Assume problem (\ref{GMP}) is feasible and the following Slater-type condition holds:
\begin{align*}\label{eqCQ}
\text{there exist scalars } z_0,z_1,\ldots,z_N  \in \bbR \text{ such that } \sum_{i=0}^{N} z_i f_i(x)>0 \text{ for all } x \in K.
\end{align*}
Then (\ref{GMP}) has an optimal solution $\mu$, which can be chosen to be \emph{finite atomic}, i.e., $\mu = \sum_{j \in J} c_j \delta_{x^{(j)}}$ for some finite index set $J$, scalars $c_j\geq 0$, and points $x^{(j)} \in K$. 
If, in addition, $\mcal{M}({\bf g})$ is Archimedean, i.e., $R-\sum_{i=1}^m x_i^2\in \mcal{M}({\bf g})$ for some scalar $R>0$, then  we have $\lim_{t \to \infty}\xi_t= \xi_{\infty}=\val$.
\end{theorem}
Hence, the link between the GMP (\ref{GMP}) and the hierarchy (\ref{hier}) is established.

Earlier, we said that, for each $t$, problem (\ref{hier}) is an SDP.
This fact may become apparent after the following characterizations.
Firstly, observe that a polynomial $\sigma \in \bbR[x]_{2t}$ is a sum of squares, i.e. $\sigma \in \Sigma$, if and only if there exists some matrix $M_{\sigma} \succeq 0$ such that $\sigma(x) = [x]_t^TM_{\sigma}[x]_t$. Having this in mind, we see that $L \geq 0$ on $\Sigma_{2t}$ is equivalent to saying $L([x]_t[x]_t^T) \succeq 0$, because $L(\sigma) = L([x]_t^TM_{\sigma}[x]_t) = \inp{L([x]_t[x]_t^T),M_{\sigma}}$, and using the fact that the PSD cone is self-dual. Similarly, for any $\sigma \in \Sigma_{2t}$ and $j \in [p]$ we have $L(g_j\sigma) = \inp{L(g_j(x)[x]_t[x]_t^T),M_{\sigma}}$. Thus
$L \geq 0$ on $\mcal{M}({ \bf g })_{2t}$ can be equivalently  characterized by the PSD constraints: 
\begin{equation}  \label{quad_mod_SDP}
L(g_j[x]_{t-d_{g_j}}[x]_{t-d_{g_j}}^T) \succeq 0,
\end{equation}  
for $j = 0,1,...,p$, where $d_{g_j}:= \lceil \deg(g_j)/2 \rceil$.
Secondly, the ideal constraints $L = 0$ on $\mcal{I}({ \bf h })_{2t}$ can be encoded as follows:
\begin{equation} \label{ideal_SDP}
L(h_j[x]_{2t-\deg(h_j)})  = 0,
\end{equation}
for each  $j \in [q]$, where the vector equality should be understood entry-wise.

\begin{example}{Exercise 2} 
Prove the above two claimed equivalences.
\end{example}
Using \cref{quad_mod_SDP} and \cref{ideal_SDP} we can reformulate problem (\ref{hier}) as
\begin{equation}  \label{hieralt}
\begin{split} 
\xi_t =  \min \{ L(f) :& L \in \bbR[x]^*_{2t},\\
			 & L(f_i) = a_i ~~(i \in [N]),\\
			 & L(g_j[x]_{t-d_{g_j}}[x]_{t-d_{g_j}}^T) \succeq 0 ~(j = 0,1,...,p), \\
			 & L(h_j[x]_{2t-\deg(h_j)})  = 0  ~(j \in [q])\}.
\end{split}
\end{equation}

For fixed level $t$, the program (\ref{hieralt}) is an SDP of size polynomial in $m$. It is known that SDPs are efficiently solvable under some technical conditions, see \cite{doi:10.1137/1.9781611970791}. However, computing $\val$ remains inefficient because the matrices describing (\ref{hieralt}) could be of size $\max_{j \in [p]}\{{m+t-d_{g_j} \choose t-d_{g_j}}\}$, and hence soon grow beyond what most currently available hardware can store in memory. In our experience the level $t$ of the hierarchy is often quite small (1 to 5) for practical examples.

\paragraph{\textbf{Finite convergence and recovering optimal solutions}}
Thus far in this section, we have seen how to get successive approximations of $\val$. We saw in the preceding \cref{sec:2:1} how GMPs relate to POPs, and soon in \cref{sec:2.3} we will see how GMPs relate to the CP-rank of a matrix.
However, what was not shown is whether we can also recover an optimizer $x^{\mathrm{min}}$. On another note, we said that the hierarchy quickly exceeds hardware capacity as the level increase, so it would be helpful if we had finite convergence, i.e., $\xi_t = \xi_\infty$ for some $t < \infty$. 
It turns out that there is a condition under which we solve both of these shortcomings mentioned above: finite convergence and the possibility of recovering an optimizer. Simply put, if a solution $L$ to the hierarchy $\xi_t$ at level $t$ satisfies the \emph{flatness condition} (\ref{eqflat}), then the bound at that level is exact, i.e., $\xi_t = \xi_\infty$, and there is a way to extract a finite atomic solution to the GMP in \cref{GMP}. We formally state the classical theorem due to Curto and Fialkow \cite{Curto1996SolutionOT,Curto2000TheTC}, upon which we base these claims. Note that the original formulation by Curto and Fialkow was not in the context of GMPs. We refer the reader again to \cite{doi:10.1142/p665, Klerk2019ASO} for a more cohesive view. 

\begin{theorem} \textbf{(Flatness theorem)} \cite{Curto1996SolutionOT,Curto2000TheTC} \label{CurtoFialkowr_flat_ext_atom} 
Consider the set $K$ from (\ref{semi_alg}) and define $d_K:=\max\{1, \lceil\deg(g_j)/2\rceil: j\in [p]\}$.
Let $t\in \bbN$ such that  $2t\ge \max\{\deg(f_i): 0\le i\le N\}$ and  $t\ge d_K$.
Assume  $L \in \bbR[x]_{2t}^*$ is an optimal solution to the program (\ref{hieralt}) defining the parameter $\xi_t$ and it satisfies the following {\em flatness condition}:
\begin{equation}\label{eqflat}
\rank\ L([x]_s[x]_s^T) =\rank\ L([x]_{s-d_K}[x]_{s-d_K}^T) \text{ for some integer } s  \text{ such that } d_K\le s\le t.
\end{equation}
Then equality $\xi_t=\val$ holds and problem (\ref{GMP}) has an optimal solution $\mu$ which is finite atomic and supported on $r:= \rank\ L([x]_{s-d_K}[x]_{s-d_K}^T)$ points in $K$.
\end{theorem}
For the details on how to extract the atoms of the optimal measure when the flatness condition of Theorem \ref{CurtoFialkowr_flat_ext_atom} holds, we refer to \cite{Henrion2005, Laurent2009}. 

In the context of factorization, we will soon see that the convergence described in Theorems \ref{theoconvGMP} and \ref{CurtoFialkowr_flat_ext_atom} is not towards the rank but instead another closely related convex parameter. Furthermore, for the cases of nonnegative- and CP factorization, the atoms we recover are exactly the columns of the factorization matrices. We now proceed to apply the above techniques to the task of approximating the completely positive rank.

\subsection{Constructing a hierarchy of lower bounds for CP-rank}
\label{sec:2.3}
With the moment hierarchy machinery in place, we return our attention to factorization ranks. In particular, we will construct a hierarchy of lower bounds for the CP-rank. It should be clear to the reader how to extend the contents of this section to the nonnegative and separable ranks. As for the other ranks discussed in \cref{sec:1}, some technicalities will be required, which we omit for brevity and simply provide references where appropriate.

Begin by recalling the definition of the CP-rank from \cref{cprankdef}, and note that we can equivalently express it as follows:
\begin{equation} \label{nnrankdefatom}
\cprank (M) = \min \{r \in \mathbb{N} : M = \sum_{\ell=1}^r a_\ell a_\ell^T \text{ for some } a_1,a_2,...,a_r \in \bbR^{m}_+  \}.
\end{equation}
The above is sometimes called an \emph{atomic} formulation because the factors $a_\ell$ (the columns of $A$ in the \cref{cprankdef} formulation) can be thought of as the atoms of the factorization. Assuming we know that $M$ is CP, we believe solving the optimization problem (\ref{nnrankdefatom}) is still hard, though, to the best of our knowledge, there is no proof of this claim. It is natural to ask if relaxing some constraints yields an easier problem. In this vein, Fawzi and Parrilo \cite{fawziSelfscaledBoundsAtomic2016a} introduced a natural ``convexification" of the CP-rank:
\begin{equation} \label{tauCP}
\tcp(M) := \inf \big{\{}  \lambda : \frac{1}{\lambda} M \in \conv\{xx^T : x \in \bbR^{m}_+,~ M-xx^T \succeq 0,~ M \geq xx^T \}  \big{\}}.
\end{equation}
A similar parameter can be defined for the NN-rank \cite{fawziSelfscaledBoundsAtomic2016a} and the separable rank \cite{Gribling2022BoundingTS}.
\begin{example}{Exercise 3} 
Prove that $\tcp(M)$ is a lower bound for $\cprank (M)$.
\end{example}

Because $\tcp$ is a convex relaxation of the combinatorial parameter $\cprank$, it is possibly strictly worse, i.e., $\tcp(M) < \cprank(M)$ for some $M$. Furthermore, $\tcp(M)$ does not appear any easier to compute than $\cprank (M)$, in part because we do not have an efficient characterization of the convex hull described in \cref{tauCP}. However, not all is lost, as $\tcp(M)$ can be reformulated as a GMP,
\begin{equation} \label{GMP_NN}
\tcp(M) = \inf_{\mu\in \mcal{M}(K^M)}\Big\{\int_{K^M}1 d\mu: \int_{K^M} x_ix_jd\mu = M_{ij}\ (i,j\in [m]),\Big\},
\end{equation}
where
\begin{equation}  \label{KM}
\begin{split} 
K^M := \{x \in \bbR^{m} :& \sqrt{M_{ii}}x_i - x_i^2 \geq 0 ~~ (i \in [m]), \\
& M_{ij} - x_i x_j \geq 0 ~~ (i \neq j \in [m]), \\ 
& M - xx^T \succeq 0\}.
\end{split}
\end{equation}
For a small proof using Theorem \ref{theoconvGMP} see Lemma 2 of \cite{Korda2022ExploitingII}. The idea of using a GMP to model CP matrices was already explored in the work of Nie \cite{10.1007/s10208-014-9225-9}.

The reader may wonder why the constraints $M_{ij} - x_i x_j \geq 0 $ and $ \sqrt{M_{ii}}x_i - x_i^2 \geq 0$ are preferred here over the equivalent and more intuitive constraints:  $x_i \geq 0$ and $M_{ij} - x_i x_j \geq 0 $ (for all $i,j\in [m]$).
This is because the former gives, for finite $t$, a larger truncated quadratic module, which in turn gives better bounds for the finite levels of the hierarchy. 
Both options are, of course, equivalent in the limit as $t$ goes to infinity, see \cite{Gribling2019LowerBO}.

Note that the last constraint, $ M - xx^T \succeq 0$, is a polynomial matrix constraint.
The idea behind this constraint is to encode that any CP factor $aa^T$ of $M$ is PSD less than $M$, i.e., $M-aa^T \succeq 0$. 
We could equivalently have asked that $f_v(x):= v^T(M - xx^T)v \geq 0$ for all $v \in \bbR^n$, or that the minors of $M - xx^T$ be nonnegative. However the matrix formulation is computationally easier to implement as we will see below. Several characterizations and supplementary references are considered for this constraint in \cite{Gribling2022BoundingTS}.

Now we simply apply the techniques of \cref{sec:2:1} to construct a hierarchy of lower bounds for the above GMP (\ref{GMP_NN}) to obtain the following parameter for any $t\in \bbN \cup \{\infty\}$:

\begin{equation}  \label{CPhier}
\begin{split} 
\xicp_t(M):=  \min \{ L(1) :& L \in \bbR[x]^*_{2t},\\
			 & L(xx^T) = M, \\
			 & L([x]_t[x]_t^T ) \succeq 0, \\
			 & L((\sqrt{M_{ii}}x_i - x_i^2)[x]_{t-1}[x]_{t-1}^T) \succeq 0  ~(i \in [m]), \\
			 & L((M_{ij} - x_ix_j)[x]_{t-1}[x]_{t-1}^T) \succeq 0 ~(i \neq j \in [m]), \\
			 & L( (M - xx^T) \otimes [x]_{t-1}[x]_{t-1}^T) \succeq 0 \}.
\end{split}
\end{equation}
The basic idea for the construction of this hierarchy comes initially from \cite{Gribling2019LowerBO}. The last constraint was added later in \cite{Gribling2022BoundingTS}.

Using Theorem \ref{theoconvGMP} we have the following chain of inequalities:
$$
\xicp_1(M) \leq \xicp_2(M) \leq ... \leq \xicp_\infty(M) = \tcp(M) \leq \cprank(M). 
$$

Let us get some intuition for why this hierarchy works. Consider a CP factorization $a_1,...,a_r \in \bbR^m_+$ of $M$ with $r:=\cprank(M)$. Define for each $i \in [r]$ the following \emph{evaluation linear functional} $L_{a_i}$ that maps a polynomial $f(x)$ to its evaluation at the point $a_i$, i.e.,
$$
L_{a_i}: \bbR[x] \ni f(x) \mapsto f(a_i) \in \bbR  ~~ (i \in [r]).
$$
Then $\tilde{L} := \sum_{i \in [r]} L_{a_i}$ is feasible for the problem (\ref{CPhier}). Moreover, $\tilde{L}(1) = r = \cprank(M)$.

\begin{example}{Exercise 4}
Show that $\tilde{L}$ satisfies each of the constraints of problem (\ref{CPhier}).
\end{example}
One can think of the constraints in (\ref{CPhier}) as filters excluding solutions that are dissimilar to $\tilde{L}$.
Of course, $\tilde{L}$ is only feasible and not necessarily optimal.
Hence $\xi_t(M) \leq \cprank(M)$ for every $t$, and in practice, the inequality is often strict. 
The finite atomic measure $\tilde\mu:= \sum_{i=1}^r \delta_{a_i}$ supported on the atoms $a_1,...,a_r$ is a representing measure of $\tilde{L}$. \\

With a hierarchy constructed, we can now compute some examples.

\subsection{A note on computing hierarchies of SDPs}
We said before that, for a fixed $t\in \bbN$, problem (\ref{CPhier}) is an SDP, and we claimed that it could be computed efficiently. We now give some tips in implementing these problems using freely available software. Theory is often a poor substitute for hands-on experience when it comes to implementing code. Therefore, this small subsection is simply  a snapshot of the quickly changing available tools. The reader is encouraged to play around with these tools should he/she seek a deeper understanding. 
At the end, we list a table of results so that the reader can get a feel for the power of these approximation hierarchies.

Disclaimer, the procedure we describe here is based on the author's experience and preferences and should not be seen as the only way to compute hierarchies.

 The core idea is to work inside the programming language Julia \cite{doi:10.1137/141000671}, within which there is a package called JuMP \cite{Dunning2017JuMPAM} specially designed as a high-level interface between several commonly used solvers and Julia. In particular, JuMP can interface with MOSEK  \cite{Andersen2000}, a powerful commercial \emph{interior-point solver}. Though MOSEK requires a license to run, academic licenses are available free of charge at the time of writing this. To summarize, one installs Julia, imports JuMP, uses the JuMP syntax to formulate the desired SDP as a JuMP-model, and then one passes off the JuMP-model to MOSEK to be solved. In broad strokes, this was the procedure followed in \cite{Gribling2019LowerBO,Gribling2022BoundingTS,Korda2022ExploitingII} to compute bounds for the NN-, CP-, and separable ranks. Some code is available as a package \footnote{See the code repository: https://github.com/JAndriesJ/ju-cp-rank}.

\paragraph{\textbf{Some numerical results for CP-rank}}
In \cite{Gribling2022BoundingTS}, the hierarchy (\ref{CPhier}) was tested on several matrices with high CP-rank. See \cite{BOMZE2014208} for the construction and definitions of these matrices. We already saw $\wtl M_7$ above in \cref{M7}. We now list the bound at level $t=3$ of the hierarchy (\ref{CPhier}) for some of the other matrices in \cite{BOMZE2014208}.

\begin{table}[!htbp]\label{CP_bounds}
	\centering
	\caption{Bounds for completely positive rank at level $t=3$.}
	\begin{tabular}{|c c c c  c c  |} 
		\toprule
		$M$ &  $\rank(M)$ &$m$& $\lfloor \frac{m^2}{4}\rfloor$  & $\xicp_{3}(M)$ & $\cprank(M)$ \\ [0.5ex] 
		\hline\hline
		$M_7$ & 7&7&12&{11.4} & 14    \\ 
		\hline
		$\wtl M_7$  &7&7&12&10.5 &  14    \\  
		\hline 
		$\wtl M_8$   &8&8&16&{14.5} & 18     \\ 
		\hline  
		$\wtl M_9$   &9&9&20&{18.4}&  26    \\  
		\bottomrule
	\end{tabular}
\end{table}
Level $t=3$ was the highest that could be computed on the available hardware.

\section{Exploiting sparsity}
\label{sec:3}
In this penultimate section, we explore \emph{sparsity}, by which we mean the exploitation of zeros in the matrix $M$ to obtain better bounds or faster computations. In particular, we will explore a special kind of sparsity called \emph{ideal sparsity}, as defined in \cite{Korda2022ExploitingII}. Recall in problem (\ref{hier}) that the ideal constraint in the SDP forces the measure to vanish on certain polynomials. In this section, we show how the zero entries in a matrix lead the measure vanish on particular monomials.
Using this, we can replace the original measure by multiple measures, each with a smaller support than the original one. The motivation behind this divide-and-conquer tactic is that the measures with smaller support lead to SDPs with smaller matrices and hence are easier to compute. Surprisingly, the sparse hierarchy, which we will define in (\ref{sparhier}), is also stronger than its dense analog (\ref{hier}), in contradistinction to other sparsity techniques where one often sacrifices the quality of the bounds in favor of computational benefits.

We first begin with a general introduction to ideal sparsity for the GMP setting in \cref{sec:31}. With the basic idea established, we apply ideal sparsity to CP-rank in \cref{sec3.2} and construct a sparse analog to the hierarchy (\ref{CPhier}). Finally, we conclude this section with some results in \cref{sec3.3} demonstrating the benefits of this sparse hierarchy over its dense analog.

\subsection{An abbreviated introduction to ideal sparsity}
\label{sec:31}
Let $V := [m]$, $E \subseteq \{\{i,j\} \in V \times V : i \neq j \}$, and let $\overline{E}:=  \{\{i,j\} \in V \times V : \{i,j\} \notin E,~i\neq j \}$ be its complement.
Suppose now that the semi-algebraic set from \cref{semi_alg} is defined as follows
\begin{equation*} 
K_E := \{ x \in \bbR^{m} : g_i(x) \geq 0 ~(i \in [p]) ~,~x_ix_j = 0 ~(\{i,j\}\in \overline{E}  ) \}.
\end{equation*}

By definition, the ideal in \cref{ideal} becomes
\begin{equation}\label{E_ideal}
\mcal{I}_{E,2t} :=  \Big\{\sum_{\{i,j\}\in \overline{E}} \gamma_{ij} x_ix_j: \gamma_{ij}\in \bbR[x]_{2t-2}\Big\}\subseteq\bbR[x]_{2t}.
\end{equation} 
Observe that we have $K_E \subseteq \mcal{I}_{E}$. We plan to partition $K_E$ in a particular way to eliminate the need for ideal constraints in the subsequent levels of the hierarchy.

Consider the undirected graph $G=(V,E)$ and denote its maximal cliques by $V_1,\ldots,V_s$. For each $k\in [s]$ define  the following subset of $K_E$ :
\begin{equation}\label{eqwhKk}
\wh {K_k} :=\{x \in K: \supp(x)\subseteq V_k \}\subseteq K\subseteq \bbR^m.
\end{equation}
Here, $\supp(x)=\{i\in [m]: x_i \neq 0\}$ denotes the {\em support} of $x\in \bbR^m$.
Observe that the sets $\wh{K_1},\ldots,\wh{K_s}$ cover the set $K_E$:
\begin{equation}\label{eqKkcover}
K_E=\wh{K_1}\cup \ldots \cup \wh{K_s}. 
\end{equation}
\noindent
It is an easy exercise to see that if $x\in K_E$, then its support $\supp(x)$ is a clique of the graph $G$, and thus it is contained in a maximal clique $V_k$,  so that  $x\in \wh{K_k}$, for some $k\in [s]$. Now  define $K_k\subseteq \bbR^{|V_k|}$ to be the projection of $\wh{K_k}$ onto the subspace indexed by $V_k$:
 \begin{align}\label{eqKk}
K_k:=\{ y \in \bbR^{|V_k|}: (y,0_{V\setminus V_k})\in \wh{K_k}\} \subseteq \bbR^{|V_k|}.
\end{align}
We use the notation $(y,0_{V\setminus V_k})$ to denote the vector of $\bbR^n$ obtained from $y\in \bbR^{|V_k|}$ by padding it with zeros at all entries indexed by $V\setminus V_k$.
For an $n$-variate function $f:\bbR^{|V|}\to\bbR$ and a subset $U\subseteq V$, we let $f_{|U}:\bbR^{|U|}\to\bbR$ denote the function in  the variables $x(U)=\{x_i: i\in U\}$, which is obtained from $f$ by setting to zero all the variables $x_i$ indexed by $i\in V\setminus U$. That is, $f_{|V_k}(y)=f(y,0_{V\setminus V_k})$ for $y\in\bbR^{|V_k|}$.

We may now define the following sparse analog of (\ref{GMP}):

\begin{equation}\label{SparGMP}
\begin{split} 
\valsp:= \inf_{\mu_k\in \meas(K_k), k\in [s]} \Big\{\sum_{k=1}^s \int  {f_0}_{|V_k}d\mu_k~:&\sum_{k=1}^s \int {f_i}_{|V_k} d\mu_k=a_i \ (i\in [N])\Big\}.
\end{split} 
\end{equation}

\begin{proposition}\label{propequiv} \cite{Korda2022ExploitingII}
Problems (\ref{GMP}) (using $K_E$) and (\ref{SparGMP}) are equivalent, i.e., their optimum values are equal: $\val=\valsp$.
\end{proposition}

Based on the reformulation (\ref{SparGMP}) we can define the following {\em ideal-sparse} moment relaxation for problem (\ref{GMP}): for any integer $t \in \bbN \cup \{\infty \}$
\begin{equation}\label{eqxisptb}
\begin{array}{ll}
\xisp_t:=\inf \Big \{\sum_{k=1}^s L_k({f_0}_{|V_k}): & L_k \in \bbR[x(V_k)]_{2t}^* \ (k\in [s]),\\
													 & \sum_{k=1}^s L_k({f_i}_{|V_k})=a_i\ (i\in [N]),\\
													 & L_k\ge 0\text{ on } \mcal{M}({\bf g}_{|V_k})_{2t} \ (k\in [s])\Big\},
\end{array}
\end{equation}
where ${ \bf g }_{|V_k}:= (g_{0|V_k},g_{1|V_k},g_{2|V_k},...,g_{p|V_k})$. Note that the ideal constraint are entirely captured by the fact that none of the measures $\mu_k$ support any elements of the ideal.

With two hierarchies converging to the same value, the obvious question is whether one converges faster. Surprisingly, the bounds for the sparse hierarchy (\ref{eqxisptb}) are at least as good as for the dense hierarchy (\ref{hier}).

\begin{theorem}\label{lem:dens<spar} \cite{Korda2022ExploitingII}
 For any integer $t \in \bbN \cup \{\infty \}$, we have $\xi_t\le \xisp_t \le \val$. If, in addition $\mcal{M}(\bf{g})$ is Archimedian and the condition in Theorem \ref{theoconvGMP} holds, then $\lim_{t\to\infty} \xisp_t = \val$.
\end{theorem}

The advantage of (\ref{eqxisptb}) over (\ref{hier}) is twofold. Firstly, the sparse bounds are at least as good as the dense bounds. Secondly, there is potential for computation speed-up since each set $V_k$ can be much smaller than the whole set $V$. This holds despite there now being more variables and constraints overall. However, speed-up fails in cases where there are exponentially many (in $n$) maximal cliques, like when $G$ is a complete graph with a perfect matching deleted.

Observe that chordality need not be assumed on the cliques. However, we are required to find all maximal cliques. For an arbitrary graph, this could be difficult, but in the setting of factorization ranks, the graphs are often small, with around $5$ to $15$ vertices. Hence, one can compute the maximal cliques using algorithms like the one described in \cite{10.1145/362342.362367}.

\subsection{Ideal sparsity in approximating CP-rank}
\label{sec3.2}
Return now to the completely positive rank. The rather abstractly defined ideal constraint in \cref{sec:31} will emerge naturally from the zeros in a matrix. Consider a CP matrix $M \in \mcal{S}_+^m$, assume $M_{ii}>0$ for all $i\in [m]$. If $M$ is a  CP matrix with $M_{ii}=0$, then its $i^\mathrm{th}$ row and column are identically zero, and thus it can be removed without changing the CP-rank.
Define the  {\em support graph} $G_M:=(V,E_M)$ of $M$, with edge-set and non-edge-set respectively defined by:  
\begin{equation*}
E_M:=\{\{i,j\}: M_{ij} \ne 0,\, i,j\in V,\, i\ne j\},\ \overline{E}_M:=\{\{i,j\}: M_{ij}=0,\, i,j\in V,\, i\ne j\}.
\end{equation*}
If $G_M$ is not connected, then $M$ can be block-diagonalized using row and column permutations. It is immediately apparent that the CP-rank of a block diagonal matrix is the sum of the CP-ranks of its blocks. So we may assume that $G_M$ is connected.
Now we can modify the semi-algebraic set from \cref{KM} to read as follows
\begin{equation}\label{eqsetKA}
\begin{array}{lll}
K_M^{\mathrm{isp}}:=\{x\in \bbR^m:   &\sqrt{M_{ii}}x_i-x_i^2\ge 0 \ (i\in [m]), \\
 & \ M_{ij}-x_ix_j\ge 0 \ (\{i,j\}\in E_M),\\
 & \ x_ix_j=0 \ (\{i,j\}\in  \overline{E}_M),\\
 &  M-xx^T\succeq 0 \}.
\end{array}
\end{equation}
We have not introduced any new information. We have just explicitly encoded the fact that $M_{ij}=0$ and $M_{ij}-x_ix_j\ge 0 $ imply $x_ix_j=0$ (because $x \geq 0$). In this form we can apply the results from \cref{sec:2.2,sec:2.3,sec:31} to define the following new hierarchy,
 
\begin{equation} \label{sparhier}
\begin{split}
\xicpsp_t(M)= \min\Big\{\sum_{k=1}^s L_k(1): & \ L_k\in \bbR[x(V_k)]^*_{2t}\ (k\in [s]),\\
&   \sum_{k\in [s]: i,j\in V_k}  L_k(x_ix_j)=M_{ij}\ (i,j\in V), \\
& L_k([x(V_k)]_t[x(V_k)]_t^T)\succeq 0 \ (k\in [s]), \\
&L_k((\sqrt{M_{ii}}x_i-x_i^2)[x(V_k)]_{t-1}[x(V_k)]_{t-1}^T)\succeq 0 \  (i\in V_k,\ k\in [s]), \\
&L_k((M_{ij}-x_ix_j)[x(V_k)]_{t-1}[x(V_k)]_{t-1}^T)\succeq 0 \ ( i\ne j\in V_k,\ k\in [s]), \\
& L_k( {(M- x x^T)} \otimes   [x(V_k)]_{t-1}[x(V_k)]_{t-1}^T)\succeq 0,\ ( k\in [s]).  
\end{split}
\end{equation} 
As a direct consequence of Theorem \ref{lem:dens<spar} we have the following relation:
$$
\xicp_t(M) \leq \xicpsp_t(M) \leq \tcp(M).
$$
 
Problem (\ref{sparhier}) looks cumbersome. However, it is ultimately just problem (\ref{CPhier}) with the single functional replaced by multiple functionals, each with support tailored to exclude polynomials in the ideal $\mcal{I}_{E_M}$.

Observe that, if, in problem (\ref{sparhier}), we replace the matrix $M-xx^T$ by its principal submatrix indexed by $V_k$, then  one also gets a lower bound on $\tau_{\mathrm{cp}}(M)$, at most $\xicpsp_t(M)$, but potentially cheaper to compute. We let $\xicpwsp_t(M)$ denote the parameter obtained in this way, by replacing in the definition of $\xicpsp_t(M)$ the last constraint by
\begin{align}\label{eqcpmatwsp}
L_k( {(M[V_k]- x(V_k)x(V_k)^T)} \otimes   [x(V_k)]_{t-1}[x(V_k)]_{t-1}^T)\succeq 0\ \text{ for } k\in [s],
\end{align}
so that we have 
 $$\xicpwsp_t(M)\le \xicpsp_t(M).$$

\begin{example}{An example of maximal cliques in the support graph of a matrix} 
To get some intuition into what the maximal cliques look like in the CP factorization setting,
 consider the following matrix and its associated support graph in \cref{cliqyy}.
\begin{figure}[ht]
        \centering
        \includegraphics[scale=0.28]{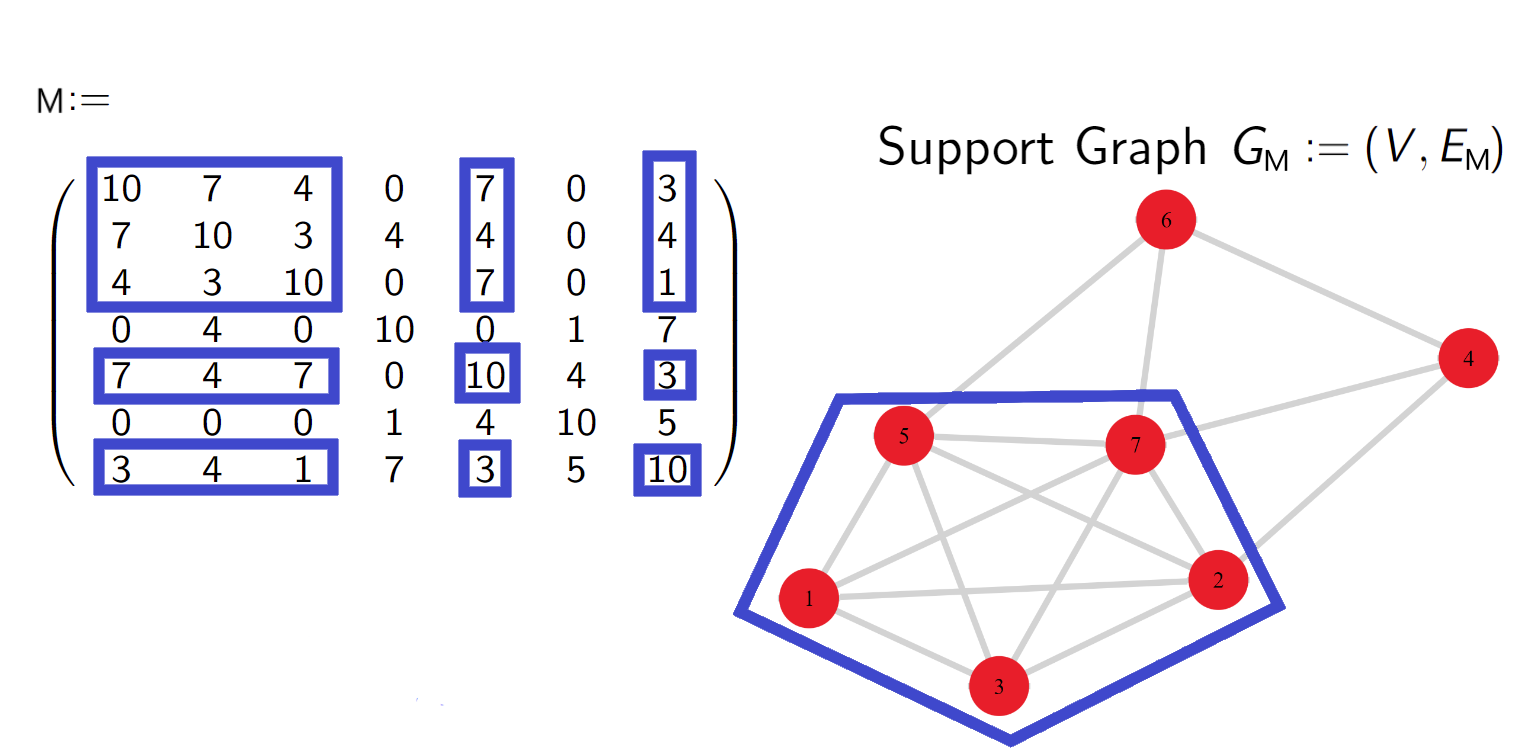}
        \caption{Example of a matrix and its support graph. This example has non-edges:
$\overline{E}_M := \{ \{1,4\}, \{1,6\}, \{2,6\}, \{3,4\}, \{3,6\}, \{4,5\} \}$, 
and maximal cliques:
{\color{blue} $V_1:= \{1,2,3,5,7\}$},  $V_2:= \{2,4,7\}$, $V_3:= \{5,6,7\}$,  $V_4:= \{4,6,7\}$. 
Hence, if $M$ is CP, then its factors can only be supported by one of these four cliques.
}
\label{cliqyy}
\end{figure}
\end{example}

\subsection{Advantages of the sparse hierarchy}
\label{sec3.3}
In this subsection, we compare the dense and sparse hierarchies for approximating the CP-rank. The comparison is first made in terms of bounds and then in terms of computational speed-up.

\paragraph{\textbf{Better bounds}}
We now demonstrate some advantages of the sparse hierarchy (\ref{sparhier}) above its dense counterpart in (\ref{CPhier}). To this end consider one of the matrices from \cite{doi:10.1137/140973207}, namely,
{\small
\[
  \widehat{M} =
  \left[ {\begin{array}{ccc ccc ccc ccc}
    91& 0&  0&  0&  19& 24& 24& 24& 19& 24& 24& 24 \\
	0&  42& 0&  0&  24& 6&  6&  6&  24& 6& 6& 6 \\
	0&  0&  42& 0&  24& 6&  6&  6&  24& 6& 6& 6 \\
	0&  0&  0&  42& 24& 6&  6&  6&  24& 6& 6& 6 \\
	19& 24& 24& 24& 91& 0&  0&  0&  19& 24& 24& 24 \\
	24& 6&  6&  6&  0&  42& 0&  0&  24& 6& 6& 6 \\
	24& 6&  6&  6&  0&  0&  42& 0&  24& 6& 6& 6 \\
	24& 6&  6&  6&  0&  0&  0&  42& 24& 6& 6& 6 \\
	19& 24& 24& 24& 19& 24& 24& 24& 91& 0& 0& 0 \\
	24& 6&  6&  6&  24& 6&  6&   6& 0& 42& 0& 0 \\
	24& 6&  6&  6&  24& 6&  6&   6& 0& 0& 42& 0 \\
	24& 6&  6&  6&  24& 6&  6&   6& 0& 0& 0& 42 \\
  \end{array} } \right].
\]}

For this matrix we know that $\cprank(\widehat{M}) = 37$. At the first level $t=1$, we have $\xicpsp_1(\widehat{M}) = 29.66$ while the dense hierarchy gives $\xicp_1(\widehat{M}) = 4.85$. Going to higher levels does improve the dense bound to $\xicp_2(\widehat{M}) =29.66$, but the sparse bound does not seem to change.

In \cite{Korda2022ExploitingII}, it was shown that the separation between $\xicpsp_1(M)$ and $\xicp_1(M)$ could be made arbitrarily big by taking matrices $M$ of the form:
$$
M=\left(\begin{matrix} (m+1)I_m & J_m\cr J_m & (m+1)I_m\end{matrix}\right)\in \mathcal{S}^{2m}
$$ 
and increasing $m$. Here $I_m$ is the identity matrix, and $J_m$ is the all-ones matrix.

 This gap is motivated by the sparse hierarchy incorporating certain structural information that the dense hierarchy ignores. To understand what we mean, consider first the \emph{edge clique-cover number} $c(G)$ of the graph $G$ defined to be the minimal number of cliques needed to cover all edges of $G$. In Lemma 13 of \cite{Korda2022ExploitingII} it is shown that
 $$
 \xicpwsp_1(M)\ge c_{\text{\rm frac }}(G_M),
 $$ 
 where $c_{\text{\rm frac }}(G)$ is the \emph{fractional edge clique-cover number}, the natural linear relaxation of $c(G)$:
 
\begin{equation*}
c_{\text{\rm frac }}(G):=\min\Big\{\sum_{k=1}^s x_k: \sum_{k: \{i,j\}\subseteq V_k} x_k\ge 1\ \text{ for } \{i,j\}\in E\Big\}.
\end{equation*}

\begin{example}{Exercise 5}
Compute the fractional and usual (integer) clique-cover number of $G_{\widehat{M}}$.
\end{example}

Hence we have the following relations:

\begin{center}
    \begin{align*}
& c_{\text{\rm frac}}(G_M) \leq \xicpwsp_1(M) \leq \xicpwsp_2(M)\leq ... \leq  \xicpwsp_\infty(M) . \\
&~~~~~~~~~~~~~~~~~~~~~~~~~~~ \rotatebox{270}{$\leq$}~~~~~~~~~~~~~~~~~~~ \rotatebox{270}{$\leq$} ~~~~~~~~~~~~~~~~~~~~~~~~~~~ \rotatebox{90}{$=$} \\
&  ~~~~~~~~~~~~~~~~~~~~~  \xicpsp_1(M) \leq \xicpsp_2(M) \leq ... \leq \xicpsp_\infty(M) = \tcp(M) \leq \cprank(M)\\
&~~~~~~~~~~~~~~~~~~~~~~~~~~~ \rotatebox{90}{$\leq$}~~~~~~~~~~~~~~~~~~~ \rotatebox{90}{$\leq$} ~~~~~~~~~~~~~~~~~~~~~~~~~~~ \rotatebox{90}{$=$} \\
&~~~~~~~~~~~~~~~~~~~~~ \xicp_1(M)~~ \leq~~~ \xicp_2(M) ~~~ \leq~ ... \leq  \xicp_\infty(M) . \\
    \end{align*}
\end{center}
The weak sparse hierarchy $\xicpwsp_t$ and the dense hierarchies $\xicp_t$ are incomparable, as there are examples where one outperforms the other and vice versa.

\paragraph{\textbf{Speed-up in computation}}
We said before that the sparse hierarchy involves smaller SDPs than the dense version and, as a result, can be computed faster. 

\begin{figure}[ht]
	\centering
	\text{
	 Computation times vs. matrix size and nonzero density, level $ t=2$}
	\includegraphics[scale=0.350]{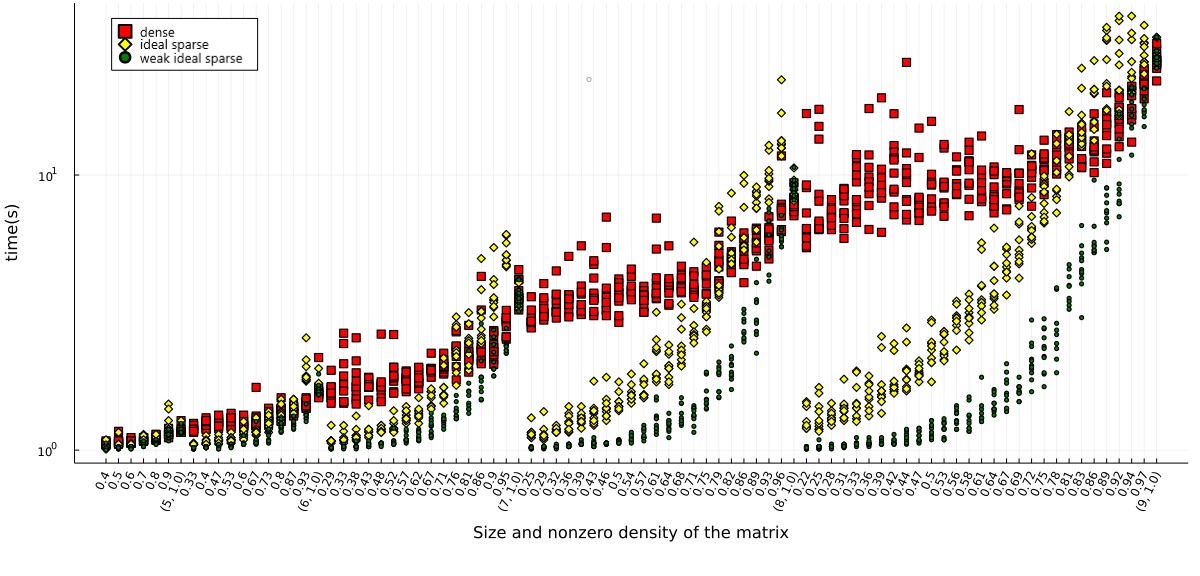}
	\caption{ Scatter plot of the computation times (in seconds) for the three hierarchies $\xicp_{2,\dag}$ (indicated by a red square), $\xicpsp_{2,\dag}$ (indicated by a yellow lozenge), $\xicpwsp_{2,\dag}$ (indicated by a green circle) against matrix size and nonzero density for 850 random matrices. The matrices are arranged in ascending size ($n=5,6,7,8,9$) and then ascending nonzero density, ranging from the minimal density needed to have a connected support graph to a fully dense matrix. For each size and nonzero density ten examples were computed to account for different support graphs.} 
	\label{scat1}
\end{figure}

To demonstrate this, the hierarchies are tested on a family of randomly generated CP-matrices ordered ascending in size and ascending in \emph{nonzero density}. The nonzero density of a matrix $M$ is the fraction of entries above the diagonal that are nonzero, hence for the identity matrix, it would be zero, and for a dense matrix with no zeros, it would be one. This parameter is crude in that it is oblivious to the structure of the support graph. Nonetheless, it suffices to show how the speed-up is related to the sparsity in the matrix.
Consider the following \cref{scat1} taken from \cite{Korda2022ExploitingII}. 

The hierarchies $\xicp_{t,\dag}$, $\xicpsp_{t,\dag}$, and $\xicpwsp_{t,\dag}$ are slight modifications of the familiar parameters $\xicp_{t}$, $\xicpsp_{t}$, and $\xicpwsp_{t}$ described already. The exact definition is avoided here because there are several technicalities to consider that will only detract from the core message, which is that the sparse hierarchy is potentially much faster when there are many zeros in the matrix.

\section{Summary}
Finally, we summarise this chapter. In \cref{sec:1}, we introduced the reader to several factorization ranks and motivated their importance with applications and links to other branches of science. After building the general tools needed, we focused on approximating the CP-rank in \cref{sec:2}. We then improved our approximation in \cref{sec:3} by including structural information about the matrix support graph before demonstrating the improvement with theoretical and numerical results. We hope to have convinced the reader of the generality and utility of polynomial optimization techniques in dealing with the difficult and pertinent problem of matrix factorization.

\paragraph{\textbf{Acknowledgements}}
We want to thank Prof. Dr. Monique Laurent for proofreading several drafts of this chapter and providing key insights when the author's knowledge was lacking. We also thank the editors for the opportunity to consolidate and share our expertise on this fascinating topic.

\bibliographystyle{abbrv} 
\bibliography{MFRchap.bib}

\end{document}